\documentclass[11pt]{article}
\usepackage{epic,latexsym,amssymb}
\usepackage{color}
\usepackage{tikz}
\usepackage{amsfonts,epsf,amsmath}

\usepackage{graphicx}
\usepackage[utf8]{inputenc}

\usepackage[english]{babel}

\usepackage{csquotes}
\usepackage{doi}

\usepackage{graphicx}
\usepackage{subfigure}
\usepackage{tikz}
\usepackage{floatrow}
\usepackage{amsthm}
\usepackage{amsfonts}
\usepackage{amsmath}
\usepackage{amssymb}
\usepackage{enumerate}
\usepackage{graphicx}
\usepackage{caption}
\usepackage{mathtools}
\usepackage{imakeidx}
\usepackage{permute}
\usepackage{array}
\usepackage{tabularx}
\usepackage{pxfonts}
\usepackage{mathtools}
\usepackage{comment}
\usepackage[noadjust]{cite}
\usepackage[sc]{mathpazo}
\usepackage[T1]{fontenc}
\usepackage{url}

\usepackage{listings}

\DeclarePairedDelimiter{\floor}{\lfloor}{\rfloor}

\newtheorem{theorem}{Theorem}[section]
\newtheorem{lemma}[theorem]{Lemma}
\newtheorem{corollary}[theorem]{Corollary}

\newtheorem{observation}[theorem]{Observation}
\newtheorem{construction}[theorem]{Construction}

\newtheorem{conjecture}[theorem]{Conjecture}

\newtheorem{problem}[theorem]{Problem}

\usepackage{color}   
\usepackage{hyperref}
\hypersetup{
    colorlinks=true, 
    linktoc=all,     
    linkcolor=blue,  
    hidelinks,
}

\newcommand{\cE}{{\cal E}}
\newcommand{\cF}{{\cal F}}

\newcommand{\cG}{{\cal G}}

\newcommand{\cH}{{\cal H}}

\let\oldenumerate\enumerate
\renewcommand{\enumerate}{
  \oldenumerate
  \setlength{\itemsep}{0pt}
  \setlength{\parskip}{0pt}
  \setlength{\parsep}{0pt}
}

\allowdisplaybreaks

\begin{document}

\title{Isolation of squares in graphs}
\author{Karl Bartolo\footnote{Email address: karl.bartolo.16@um.edu.mt} \quad \quad Peter Borg\footnote{Email address: peter.borg@um.edu.mt} \quad \quad Dayle Scicluna\footnote{Email address: dayle.scicluna.09@um.edu.mt} \\[5mm]
{\normalsize Department of Mathematics} \\
{\normalsize Faculty of Science} \\
{\normalsize University of Malta}\\
{\normalsize Malta}\\
}
\date{}

\maketitle

\begin{abstract}
Given a set $\cF$ of graphs, we call a copy of a graph in $\cF$ an $\cF$-graph. The $\cF$-isolation number of a graph $G$, denoted by $\iota(G,\cF)$, is the size of a smallest subset $D$ of the vertex set $V(G)$ such that the closed neighbourhood of $D$ intersects the vertex sets of the $\cF$-graphs contained by $G$ (equivalently, $G - N[D]$ contains no $\cF$-graph). Thus, $\iota(G,\{K_1\})$ is the domination number of $G$. The second author showed that if $\cF$ is the set of cycles and $G$ is a connected $n$-vertex graph that is not a triangle, then $\iota(G,\cF) \leq \left \lfloor \frac{n}{4} \right \rfloor$. This bound is attainable for every $n$ and solved a problem of Caro and Hansberg. A question that arises immediately is how much smaller an upper bound can be if $\cF = \{C_k\}$ for some $k \geq 3$, where $C_k$ is a cycle of length $k$. The problem is to determine the smallest real number $c_k$ (if it exists) such that for some finite set $\mathcal{E}_k$ of graphs, $\iota(G, \{C_k\}) \leq c_k |V(G)|$ for every connected graph $G$ that is not an $\mathcal{E}_k$-graph. The above-mentioned result yields $c_3 = \frac{1}{4}$ and $\mathcal{E}_3 = \{C_3\}$. The second author also showed that if $k \geq 5$ and $c_k$ exists, then $c_k \geq \frac{2}{2k + 1}$. We prove that $c_4 = \frac{1}{5}$ and determine $\mathcal{E}_4$, which consists of three $4$-vertex graphs and six $9$-vertex graphs. The $9$-vertex graphs in $\mathcal{E}_4$ were fully determined by means of a computer program. A method that has the potential of yielding similar results is introduced.
\end{abstract}


\section{Introduction}

Unless stated otherwise, we use small letters such as $x$ to denote non-negative integers or elements of sets, and capital letters such as $X$ to denote sets or graphs. The set of positive integers is denoted by $\mathbb{N}$. For $n \geq 1$, $[n]$ denotes the set $\{1, \dots, n\}$ (that is, $[n]= \{i \in \mathbb{N} \colon i \leq n\}$). We take $[0]$ to be the empty set $\emptyset$. Arbitrary sets are taken to be finite. For a set $X$, $\binom{X}{2}$ denotes the set of $2$-element subsets of $X$ (that is, $\binom{X}{2} = \{ \{x,y \} \colon x,y \in X, x \neq y \}$). We may represent a $2$-element set $\{x,y\}$ by $xy$. For standard terminology in graph theory, we refer the reader to \cite{West}. Most of the terminology used here is defined in \cite{Borg1}.

Every graph $G$ is taken to be \emph{simple}, that is, $G$ is a pair $(V(G), E(G))$ such that $V(G)$ and $E(G)$ (the vertex set and the edge set of $G$) are sets that satisfy $E(G) \subseteq \binom{V(G)}{2}$. We call $G$ an \emph{$n$-vertex graph} if $|V(G)| = n$. We call $G$ an \emph{$m$-edge graph} if $|E(G)| = m$. For any vertex $v$ of $G$, $N_G(v)$ denotes the set $\{w \in V(G) \colon vw \in E(G)\}$ of neighbours of $v$ in $G$, $N_G[v]$ denotes the closed neighbourhood $N_G(v) \cup \{ v \}$ of $v$, and $d_G(v)$ denotes the degree $|N_G (v)|$ of $v$. The minimum degree and the maximum degree of $G$ are denoted by $\delta(G)$ and $\Delta(G)$, respectively. A graph is \emph{subcubic} if its maximum degree is at most $3$. A graph is \emph{subquartic} if its maximum degree is at most $4$. For a subset $X$ of $V(G)$, $N_G[X]$ denotes the closed neighbourhood $\bigcup_{v \in X} N_G[v]$ of $X$, $G[X]$ denotes the subgraph $(X,E(G) \cap \binom{X}{2})$ of $G$ induced by $X$, and $G - X$ denotes the graph $G[V(G) \setminus X]$ obtained by deleting the vertices in $X$ from $G$. We may abbreviate $G - \{x\}$ to $G-x$. Where no confusion arises, the subscript $G$ may be omitted; for example, $N_G(v)$ may be abbreviated to $N(v)$. A \emph{component of $G$} is a maximal connected subgraph of $G$. Clearly, the components of $G$ are pairwise vertex-disjoint, and their union is $G$. 

Consider two graphs $G$ and $H$. If $G$ is a copy of $H$, then we write $G \simeq H$ and we say that $G$ is an \emph{$H$-copy}. We say that \emph{$G$ contains $H$} if $H$ is a subgraph of $G$. We say that $G$ is \emph{$H$-free} if $G$ contains no copy of $H$.  
  
For $n \geq 1$, $K_n$ and $P_n$ denote the graphs $([n], \binom{[n]}{2})$ and $([n], \{\{i,i+1\} \colon i \in [n-1]\})$, respectively. For $n \geq 3$, $C_n$ denotes the graph $([n], \{\{1,2\}, \{2,3\}, \dots, \{n-1,n\}, \{n,1\}\})$. A $K_n$-copy is called an \emph{$n$-clique} or a \emph{complete graph}, a $P_n$-copy is called an \emph{$n$-path} or simply a path, and a $C_n$-copy is called an \emph{$n$-cycle} or simply a cycle. Note that a $3$-clique is a $3$-cycle, which is also referred to as a \emph{triangle}. A $4$-cycle is also referred to as a \emph{square}. 


If $\cF$ is a set of graphs and $F$ is a copy of a graph in $\cF$, then we call $F$ an \emph{$\cF$-graph}. If $G$ is a graph and $D \subseteq V(G)$ such that $G-N[D]$ contains no $\cF$-graph, then $D$ is called an \emph{$\cF$-isolating set of $G$}. The size of a smallest $\cF$-isolating set of $G$ is denoted by $\iota(G, \cF)$ and called the \emph{$\cF$-isolation number of $G$}. If $\cF = \{F\}$, then we may replace $\cF$ in these defined terms and notation by $F$. 

The study of isolating sets was initiated in the paper \cite{CaHa17} of Caro and Hansberg. It is a natural generalization of the popular study of dominating sets \cite{C,CH,HHS,HHS2,HL,HL2}. Indeed, $D$ is a \emph{dominating set of $G$} (that is, $N[D] = V(G)$) if and only if $D$ is a $K_1$-isolating set of $G$, so the \emph{domination number of $G$}, denoted by $\gamma(G)$, is the $K_1$-isolation number of $G$. One of the earliest domination results, due to Ore \cite{Ore}, is that $\gamma(G) \leq n/2$ for any connected $n$-vertex graph $G \not\simeq K_1$ (see \cite{HHS}). While the deletion of the closed neighbourhood of a dominating set produces the graph with no vertices, the deletion of the closed neighbourhood of a $K_2$-isolating set produces a graph with no edges. Caro and Hansberg \cite{CaHa17} proved that if $G$ is a connected $n$-vertex graph, then $\iota(G, K_2) \leq n/3$ unless $G \simeq K_2$ or $G \simeq C_5$. Their seminal paper \cite{CaHa17} also featured several problems. Generalizing the bounds above, the second author of the present paper, Fenech and Kaemawichanurat \cite{BFK} showed that for any $k \geq 1$, $\iota(G, K_k) \leq n/(k+1)$ unless $G \simeq K_k$ or $k=2$ and $G \simeq C_5$ (in \cite{BFK2}, they gave a sharp upper bound in terms of the number of edges). This bound is sharp and settled a problem in \cite{CaHa17}. Domination and isolation have been particularly investigated for maximal outerplanar graphs \cite{BK, BK2, CaWa13, CaHa17, Ch75, DoHaJo16, DoHaJo17, HeKa18, LeZuZy17, Li16, MaTa96, To13, KaJi}, mostly due to connections with Chv\'{a}tal's Art Gallery Theorem \cite{Ch75}.

Let $\mathcal{C}$ denote the set of cycles. The second author \cite{Borg1} obtained a sharp upper bound on $\iota(G, \mathcal{C})$, and consequently settled another problem in \cite{CaHa17}. Before stating the result, we recall the explicit construction used for establishing an extremal case. 

\begin{construction}[\cite{Borg1}] \label{Bconstruction} \emph{For any $n, k \in \mathbb{N}$ and any $k$-vertex graph $F$, we construct a connected $n$-vertex graph $B_{n,F}$ as follows.  If $n \leq k$, then let $B_{n,F} = P_n$. If $n \geq k+1$, then let $a_{n,k} =  \big\lfloor \frac{n}{k+1} \big\rfloor$, let $b_{n,k} = n - ka_{n,k}$ (so $a_{n,k} \leq b_{n,k} \leq a_{n,k} + k$), let $F_1, \dots, F_{a_{n,k}}$ be copies of $F$ such that $P_{b_{n,k}}, F_1, \dots, F_{a_{n,k}}$ are pairwise vertex-disjoint, and let $B_{n,F}$ be the graph with $V(B_{n,F}) = [b_{n,k}] \cup \bigcup_{i=1}^{a_{n,k}} V(F_i)$ and $E(B_{n,F}) = E(P_{a_{n,k}}) \cup \{\{a_{n,k}, j\} \colon j \in [b_{n,k}] \setminus [a_{n,k}]\} \cup \bigcup_{i=1}^{a_{n,k}} (E(F_i) \cup \{\{i,v\} \colon v \in V(F_i)\})$.}
\end{construction}
\begin{theorem}[\cite{Borg1}] \label{Borgcycle}
If $G$ is a connected $n$-vertex graph that is not a triangle, then
\[\iota(G, \mathcal{C}) \leq \left \lfloor \frac{n}{4} \right \rfloor.\] 
Moreover, equality holds if $G = B_{n,C_3}$.
\end{theorem}

Problem~7.4 in \cite{CaHa17} asks for bounds on $\iota(G,\cF)$ for other interesting sets $\cF$. A question that arises immediately from Theorem~\ref{Borgcycle} is how much smaller an upper bound on $\iota(G, \cF)$ can be if instead of taking $\cF$ to be the set of cycles, we take $\cF$ to be the set of $k$-cycles. The main problem, posed in \cite{Borgrsc}, is to determine the smallest real number $c_k$ (if it exists) such that for some finite set $\mathcal{E}_k$ of graphs, $\iota(G, C_k) \leq c_k |V(G)|$ for every connected graph $G$ that is not an $\mathcal{E}_k$-graph. Also in \cite{Borgrsc}, the second author of the present paper showed that if $k \geq 5$ and $c_k$ exists, then $c_k \geq \frac{1}{k + \frac{1}{2}} = \frac{2}{2k+1}$. Recall that $C_k$ was defined for $k \geq 3$. Following \cite{Borgrsc}, we take $C_1$ and $C_2$ to be $K_1$ and $K_2$, respectively, to obtain $C_k = ( [k], \{ ij \in \binom{[k]}{2} \colon j = (i+1) \mbox{ mod$^*$ } k \} )$ for every $k \geq 1$, where mod$^*$ is the usual modulo operation except that for every $a \geq 1$ and $b \geq 0$, $ba \mbox{ mod$^*$ } a$ is $a$ instead of $0$. Ore's result yields $c_1 = \frac{1}{2}$ and $\cE_1 = \{C_1\}$, the Caro--Hansberg result yields $c_2 = \frac{1}{3}$ and $\cE_2 = \{C_2, C_5\}$, and Theorem~\ref{Borgcycle} yields $c_3 = \frac{1}{4}$ and $\mathcal{E}_3 = \{C_3\}$ (all this is also given by the result in \cite{BFK}). In this paper, we provide the solution for $k = 4$, given by Theorem~\ref{ThmC4isol}. Consequently, $c_k = \frac{1}{k+1}$ for $1 \leq k \leq 4$. As pointed out above, for any other $k$, $c_k \neq \frac{1}{k+1}$.  

Let $C_4'$ be the \emph{diamond graph} $([4],E(C_4)\cup \{\{1,3\}\})$. Let $G_1$, $G_2$, \dots, $G_6$ be the graphs whose drawings are given in Figure~\ref{FigG9}. Let $\mathcal{G}_4 = \{C_4,C_4',K_4\}$, $\mathcal{G}_9 = \{G_i \colon i \in [6]\}$ and $\mathcal{E}_4 = \mathcal{G}_4 \cup \mathcal{G}_9$. 

\begin{figure}[htb!]
\ffigbox[\textwidth]
    {
    \begin{floatrow}
    \ffigbox[\linewidth]
      {\captionof{subfigure}{$G_1$}
      \label{subfig:G1}}
      {\begin{tikzpicture}[every edge/.style={draw=black, very thick}]
          \foreach \i in {1}
            \node[shape=circle,fill,draw=black, inner sep=0pt, scale=0.5, label={north: $\i$}] (\i) at (130-40*\i:1.5cm) {\i};
        \foreach \i in {2,3,4}
            \node[shape=circle,fill,draw=black, inner sep=0pt, scale=0.5, label={east: $\i$}] (\i) at (130-40*\i:1.5cm) {\i};
            \foreach \i in {5}
            \node[shape=circle,fill,draw=black, inner sep=0pt, scale=0.5, label={south: $\i$}] (\i) at (130-40*\i:1.5cm) {\i};
            \foreach \i in {6}
            \node[shape=circle,fill,draw=black, inner sep=0pt, scale=0.5, label={south: $\i$}] (\i) at (130-40*\i:1.5cm) {\i};
            \foreach \i in {7,8,9}
            \node[shape=circle,fill,draw=black, inner sep=0pt, scale=0.5, label={west: $\i$}] (\i) at (130-40*\i:1.5cm) {\i};

            \draw (2) edge (1);
            \draw (1) edge (9);
            \draw (9) edge (8);
            \draw (8) edge (7);
            \draw (7) edge (6);
            \draw (6) edge (5);
            \draw (5) edge (4);
            \draw (4) edge (3);
            \draw (3) edge (2);
           
            \draw (2) edge (9);
            \draw (2) edge (4);
            \draw (2) edge (3);
            \draw (1) edge (8);
            \draw (1) edge (3);
            \draw (9) edge (7);  
            \draw (8) edge (6);
            \draw (7) edge (5);
            \draw (6) edge (4);  
            \draw (5) edge (3);
            
        \end{tikzpicture}}        
      
    \ffigbox[\linewidth]
      {\captionof{subfigure}{$G_2$}
      \label{subfig:G2}}
      {
      \begin{tikzpicture}[scale = 2, every edge/.style={draw=black, very thick}]
            \node[shape=circle,fill,draw=black, inner sep=0pt, scale=0.5, label={west:$5$}] (5) at (-1,0) {1};
            \node[shape=circle,fill,draw=black, inner sep=0pt, scale=0.5, label={west:$6$}] (6) at (-0.5,0.8660254038) {2};
            \node[shape=circle,fill,draw=black, inner sep=0pt, scale=0.5, label={east:$1$}] (1) at (0.5,0.8660254038) {3};
            \node[shape=circle,fill,draw=black, inner sep=0pt, scale=0.5, label={east:$2$}] (2) at (1,0) {4};
            \node[shape=circle,fill,draw=black, inner sep=0pt, scale=0.5, label={east:$3$}] (3) at (0.5,-0.8660254038) {5};
            \node[shape=circle,fill,draw=black, inner sep=0pt, scale=0.5, label={west:$4$}] (4) at (-0.5,-0.8660254038) {6} ;
            \node[shape=circle,fill,draw=black, inner sep=0pt, scale=0.5, label={west:$7$}] (7) at (0.5,0) {7};
            \node[shape=circle,fill,draw=black, inner sep=0pt, scale=0.5, label={north east:$8$}] (8) at (-0.25,-0.4330127019) {8};
            \node[shape=circle,fill,draw=black, inner sep=0pt, scale=0.5, label={south east:$9$}] (9) at (-0.25,0.4330127019) {9} ;
            
            \draw (5) edge (6);
            \draw (5) edge (8);
            \draw (5) edge (9);
            \draw (5) edge (4);
            \draw (6) edge (1);
            \draw (6) edge (7);
            \draw (6) edge (8);
            \draw (1) edge (2);
            \draw (1) edge (7);
            \draw (1) edge (9);
            \draw (2) edge (3);  
            \draw (2) edge (8);
            \draw (2) edge (9);
            \draw (3) edge (4);
            \draw (3) edge (7);
            \draw (3) edge (8);
            \draw (4) edge (7);  
            \draw (4) edge (9);
        \end{tikzpicture}
        }
    \end{floatrow}%
    
    \vspace*{15pt}

    \begin{floatrow}
    \ffigbox[\linewidth]
      {\captionof{subfigure}{$G_3$}
      \label{subfig:G3}}
      {\begin{tikzpicture}[scale =1, every edge/.style={draw=black, very thick}]
            \node[shape=circle,fill,draw=black, inner sep=0pt, scale=0.5, label={west:$5$}] (5) at (0,0.2) {5};
            \node[shape=circle,fill,draw=black, inner sep=0pt, scale=0.5, label={west:$6$}] (6) at (0,1.8) {6};
            \node[shape=circle,fill,draw=black, inner sep=0pt, scale=0.5, label={west:$7$}] (7) at (1,1) {7};
            \node[shape=circle,fill,draw=black, inner sep=0pt, scale=0.5, label={north:$8$}] (8) at (2,1) {8};
            \node[shape=circle,fill,draw=black, inner sep=0pt, scale=0.5, label={south:$4$}] (4) at (2,0) {4};
            \node[shape=circle,fill,draw=black, inner sep=0pt, scale=0.5, label={north:$1$}] (1) at (2,2) {1} ;
            \node[shape=circle,fill,draw=black, inner sep=0pt, scale=0.5, label={east:$2$}] (2) at (4,1.8) {2};
            \node[shape=circle,fill,draw=black, inner sep=0pt, scale=0.5, label={east:$3$}] (3) at (4,0.2) {3};
            \node[shape=circle,fill,draw=black, inner sep=0pt, scale=0.5, label={east:$9$}] (9) at (3,1) {9} ;
            
            \draw (5) edge (6);
            \draw (5) edge (7);
            \draw (5) edge (8);
            \draw (5) edge (4);
            \draw (6) edge (7);
            \draw (6) edge (8);
            \draw (6) edge (1);
            \draw (7) edge (1);
            \draw (7) edge (4);
            \draw (8) edge (2);
            \draw (8) edge (3);
            \draw (4) edge (3);
            \draw (4) edge (9);
            \draw (1) edge (2);  
            \draw (1) edge (9);
            \draw (2) edge (3);  
            \draw (2) edge (9);
            \draw (3) edge (9);  
        \end{tikzpicture}
}
    \ffigbox[\linewidth]
      {\captionof{subfigure}{$G_4$}
      \label{subfig:G4}}
      {\begin{tikzpicture}[scale =1, every edge/.style={draw=black, very thick}]
            \node[shape=circle,fill,draw=black, inner sep=0pt, scale=0.5, label={north:$1$}] (1) at (2,1) {1};
            \node[shape=circle,fill,draw=black, inner sep=0pt, scale=0.5, label={west:$2$}] (2) at (0,0) {2};
            \node[shape=circle,fill,draw=black, inner sep=0pt, scale=0.5, label={east:$3$}] (3) at (1,0) {3};
            \node[shape=circle,fill,draw=black, inner sep=0pt, scale=0.5, label={west:$4$}] (4) at (3,0) {4};
            \node[shape=circle,fill,draw=black, inner sep=0pt, scale=0.5, label={east:$5$}] (5) at (4,0) {5};
            \node[shape=circle,fill,draw=black, inner sep=0pt, scale=0.5, label={west:$6$}] (6) at (0,-1) {6} ;
            \node[shape=circle,fill,draw=black, inner sep=0pt, scale=0.5, label={north east:$7$}] (7) at (1,-1) {7};
            \node[shape=circle,fill,draw=black, inner sep=0pt, scale=0.5, label={north west:$8$}] (8) at (3,-1) {8};
            \node[shape=circle,fill,draw=black, inner sep=0pt, scale=0.5, label={east:$9$}] (9) at (4,-1) {9} ;
    
            \draw (1) edge (2);
            \draw (1) edge (3);
            \draw (1) edge (4);
            \draw (1) edge (5);
            \draw (2) edge (3);
            \draw (2) edge (6);
            \draw (2) edge (7);
            \draw (3) edge (6);
            \draw (3) edge (7);
            \draw (4) edge (5);  
            \draw (4) edge (8);
            \draw (4) edge (9);
            \draw (5) edge (8);  
            \draw (5) edge (9);
            \draw (6) edge (7);  
            \draw (6) edge[bend right =0.8cm] (9);
            \draw (7) edge[bend right =0.5cm] (8);
            \draw (8) edge (9);
        \end{tikzpicture}}
    \end{floatrow}
    
    \vspace*{15pt}
    
    \begin{floatrow}
    \ffigbox[\linewidth]
      {\captionof{subfigure}{$G_5$}\label{subfig:G5}}
        {\begin{tikzpicture} [scale =1, every edge/.style={draw=black, very thick}]
            \node[shape=circle,fill,draw=black, inner sep=0pt, scale=0.5, label={north:$1$}] (1) at (2,1) {1};
            \node[shape=circle,fill,draw=black, inner sep=0pt, scale=0.5, label={west:$2$}] (2) at (0,0) {2};
            \node[shape=circle,fill,draw=black, inner sep=0pt, scale=0.5, label={east:$3$}] (3) at (1,0) {3};
            \node[shape=circle,fill,draw=black, inner sep=0pt, scale=0.5, label={west:$4$}] (4) at (3,0) {4};
            \node[shape=circle,fill,draw=black, inner sep=0pt, scale=0.5, label={east:$5$}] (5) at (4,0) {5};
            \node[shape=circle,fill,draw=black, inner sep=0pt, scale=0.5, label={west:$6$}] (6) at (0,-1) {6} ;
            \node[shape=circle,fill,draw=black, inner sep=0pt, scale=0.5, label={north east:$7$}] (7) at (1,-1) {7};
            \node[shape=circle,fill,draw=black, inner sep=0pt, scale=0.5, label={north west:$8$}] (8) at (3,-1) {8};
            \node[shape=circle,fill,draw=black, inner sep=0pt, scale=0.5, label={east:$9$}] (9) at (4,-1) {9} ;
            
            \draw (1) edge (2);
            \draw (1) edge (3);
            \draw (1) edge (4);
            \draw (1) edge (5);
            \draw (2) edge (6);
            \draw (2) edge (7);
            \draw (3) edge (6);
            \draw (3) edge (7);
            \draw (4) edge (5);  
            \draw (4) edge (8);
            \draw (4) edge (9);
            \draw (5) edge (8);  
            \draw (5) edge (9);
            \draw (6) edge (7);  
            \draw (6) edge[bend right =0.8cm] (9);
            \draw (7) edge[bend right =0.5cm] (8);
            \draw (8) edge (9);
        \end{tikzpicture}}
    \ffigbox[\linewidth]
      {\captionof{subfigure}{$G_6$}\label{subfig:G6}}
      {\begin{tikzpicture} [scale =1, every edge/.style={draw=black, very thick}]
        \node[shape=circle,fill,draw=black, inner sep=0pt, scale=0.5, label={north:$1$}] (1) at (2,1) {1};
        \node[shape=circle,fill,draw=black, inner sep=0pt, scale=0.5, label={west:$2$}] (2) at (0,0) {2};
        \node[shape=circle,fill,draw=black, inner sep=0pt, scale=0.5, label={east:$3$}] (3) at (1,0) {3};
        \node[shape=circle,fill,draw=black, inner sep=0pt, scale=0.5, label={west:$4$}] (4) at (3,0) {4};
        \node[shape=circle,fill,draw=black, inner sep=0pt, scale=0.5, label={east:$5$}] (5) at (4,0) {5};
        \node[shape=circle,fill,draw=black, inner sep=0pt, scale=0.5, label={west:$6$}] (6) at (0,-1) {6} ;
        \node[shape=circle,fill,draw=black, inner sep=0pt, scale=0.5, label={north east:$7$}] (7) at (1,-1) {7};
        \node[shape=circle,fill,draw=black, inner sep=0pt, scale=0.5, label={north west:$8$}] (8) at (3,-1) {8};
        \node[shape=circle,fill,draw=black, inner sep=0pt, scale=0.5, label={east:$9$}] (9) at (4,-1) {9} ;
    
        \draw (1) edge (2);
        \draw (1) edge (3);
        \draw (1) edge (4);
        \draw (1) edge (5);
        \draw (2) edge (6);
        \draw (2) edge (7);
        \draw (3) edge (6);
        \draw (3) edge (7); 
        \draw (4) edge (8);
        \draw (4) edge (9);
        \draw (5) edge (8);  
        \draw (5) edge (9);
        \draw (6) edge (7);  
        \draw (6) edge[bend right =0.8cm] (9);
        \draw (7) edge[bend right =0.5cm] (8);
        \draw (8) edge (9);
    \end{tikzpicture}}
    \end{floatrow}
    }{\caption{The $9$-vertex graphs whose $C_4$-isolation number is $2$.}\label{FigG9}}
\end{figure}

It is straightforward that the $C_4$-isolation number of every $\mathcal{G}_4$-graph is $1$ and that the $C_4$-isolation number of every $\mathcal{G}_9$-graph is $2$. Thus, if $G$ is an $\mathcal{E}_4$-graph, then $\iota(G, C_4) > \lfloor |V(G)|/5 \rfloor$. We prove that $\mathcal{E}_4$-graphs are the only connected graphs for which this inequality holds. More precisely, in Section~\ref{proofsection}, we prove the following result. 

\begin{theorem}\label{ThmC4isol}
If $G$ is a connected $n$-vertex graph that is not an $\mathcal{E}_4$-graph, then 
\[\iota(G, C_4) \leq \left \lfloor\frac{n}{5} \right \rfloor.\]
Moreover, equality holds if $G = B_{n,C_4}$.
\end{theorem}
\noindent Yan \cite{Y} proved that $\iota(G, C_4') \leq \lfloor n/5 \rfloor$ for any connected $n$-vertex graph $G$ that is not a $\{C_4', K_4, G_1\}$-graph. This is an immediate consequence of Theorem~\ref{ThmC4isol}. 

The proof of Theorem~\ref{ThmC4isol} adopts the strategy introduced in \cite{BorgIsolConnected2021}. The challenge in proving the theorem mostly came from unexpectedly having six (non-isomorphic) $9$-vertex graphs for which the bound does not hold (a significant number of graphs with a significant number of vertices and with varying structures).  These graphs were fully determined by means of a computer program that performs an exhaustive check on $9$-vertex graphs; see the appendix. Thus, the work presented here indicates how difficult the general $C_k$-isolation problem is but may provide a way to further progress on this problem. The method introduced in this paper can be used to obtain similar results.

\section{Proof of Theorem~\ref{ThmC4isol}} \label{proofsection}

We now prove Theorem~\ref{ThmC4isol}. We often use the following two lemmas from \cite{Borg1}.

\begin{lemma}[\cite{Borg1}] \label{LemIsol1}
    If $G$ is a graph, $\cF$ is a set of graphs, $X\subseteq V(G)$, and $Y\subseteq N[X]$, then 
\[\iota(G,\cF) \leq |X| + \iota(G-Y,\cF).\]
\end{lemma}
\noindent
\textbf{Proof.} Let $D$ be an $\mathcal{F}$-isolating set of $G-Y$ of size $\iota(G-Y, \mathcal{F})$. Clearly, $\emptyset \neq V(F) \cap Y \subseteq V(F) \cap N[X]$ for each $\mathcal{F}$-graph $F$ that is a subgraph of $G$ and not a subgraph of $G-Y$. Thus, $X \cup D$ is an $\mathcal{F}$-isolating set of $G$. The result follows. \qed
\medskip

Let ${\rm C}(G)$ denote the set of components of a graph $G$.
    
\begin{lemma}[\cite{Borg1}] \label{LemIsol2}
If $G$ is a graph and $\cF$ is a set of graphs, then 
\[\iota(G,\cF) = \sum_{H \in {\rm C}(G)} \iota(H,\cF).\]
\end{lemma}

The next lemma concerns a case where no member of a subset $Y$ of $V(G)$ is a vertex of an $\cF$-graph contained by $G$, where $\cF$ is a set of cycles. 

\begin{lemma}\label{LemIsolG'}
If $G$ is a graph, $\cF$ is a set of cycles, $x \in V(G)$, $Y \subseteq V(G) \setminus \{x\}$, $N[Y] \cap V(G - Y) \subseteq \{x\}$, and $G[\{x\} \cup Y]$ contains no $\cF$-graph, then $\iota(G,\cF) = \iota(G-Y,\cF)$ and every $\mathcal{F}$-isolating set of $G-Y$ is an $\mathcal{F}$-isolating set of $G$. 
\end{lemma}

\noindent
\textbf{Proof.} Let $G' = G-Y$. Let $D$ be an $\cF$-isolating set of $G$ of size $\iota(G,\cF)$. Let $D_x = (D \setminus Y) \cup \{x\}$ if $D \cap Y \neq \emptyset$, and let $D_x = D$ otherwise. Since $N[Y] \cap V(G - Y) \subseteq \{x\}$ and $G[\{x\} \cup Y]$ contains no $\cF$-graph, $D_x$ is an $\cF$-isolating set of $G$. Since $D_x \subseteq V(G')$, $D_x$ is an $\cF$-isolating set of $G'$. Thus, we have $\iota(G',\cF) \leq |D_x| \leq |D| = \iota(G,\cF)$.

Let $D'$ be an $\cF$-isolating set of $G'$. Suppose that $G$ contains an $\cF$-graph $H$. Suppose $V(H) \nsubseteq V(G')$. Then, $V(H) \cap Y \neq \emptyset$. Since $G[\{x\} \cup Y]$ contains no $\cF$-graph, $V(H) \setminus (\{x\} \cup Y) \neq \emptyset$. Thus, $z_1z_2 \in E(H)$ for some $z_1 \in Y$ and $z_2 \in V(H) \setminus Y$. Since $N[Y] \cap V(G - Y) \subseteq \{x\}$, $z_2 = x$. We have $E(H) = \{z_1z_2, z_2z_3, \dots, z_{r-1}z_r, z_rz_1\}$ for some $r \geq 3$ and some distinct members $z_3, \dots, z_r$ of $V(G) \setminus \{z_1, z_2\}$. We have $N_H(x) = N_H(z_2) = \{z_1, z_3\}$. Suppose $z_3 \in Y$. Since $x \notin \{z_3, \dots, z_r\}$ and $N[Y] \cap V(G - Y) \subseteq \{x\}$, $z_3, \dots, z_r \in Y$. This yields $V(H) \subseteq \{x\} \cup Y$, which contradicts $V(H) \setminus (\{x\} \cup Y) \neq \emptyset$. Thus, $z_3 \notin \{x\} \cup Y$. Since $x \notin \{z_3, \dots, z_r\}$ and $N[Y] \cap V(G - Y) \subseteq \{x\}$, we obtain $z_3, \dots, z_r, z_1 \notin \{x\} \cup Y$, which contradicts $z_1 \in Y$. Therefore, $V(H) \subseteq V(G')$. Since $D'$ is an $\cF$-isolating set of $G'$, $N_{G'}[D'] \cap V(H) \neq \emptyset$. Thus, $N_{G}[D'] \cap V(H) \neq \emptyset$, and hence $G-N_G[D']$ does not contain $H$. Therefore, $D'$ is an $\cF$-isolating set of $G$. This yields $\iota(G,\cF) \leq \iota(G',\cF)$. Since we also obtained $\iota(G',\cF) \leq \iota(G,\cF)$, equality holds.\qed
\medskip

An \emph{isolated vertex of $G$} is a vertex of $G$ of degree $0$. A \emph{leaf of $G$} is a vertex of $G$ of degree $1$.

\begin{corollary}\label{LemIsolG'cor}
If $G$ is a graph, $\cF$ is a set of cycles, and $y$ is an isolated vertex of $G$ or a leaf of $G$, then $\iota(G,\cF) = \iota(G-y,\cF)$.
\end{corollary}
\noindent
\textbf{Proof.} Let $Y = \{y\}$. Suppose that $y$ is an isolated vertex of $G$. If $V(G) = Y$, then $\iota(G,\cF) = 0 = \iota(G-y,\cF)$. If $V(G) \neq Y$ and $x \in V(G) \setminus Y$, then $\iota(G,\cF) = \iota(G-y,\cF)$ by Lemma~\ref{LemIsolG'}. Now suppose that $y$ is a leaf of $G$.  Since $y$ has only one neighbour $x$, $\iota(G,\cF) = \iota(G-y,\cF)$ by Lemma~\ref{LemIsolG'}. \qed

\subsection{The case where $G$ is subcubic}\label{subcubic section}

In the proof of Theorem~\ref{ThmC4isol}, the case where $G$ is subcubic needs to be treated differently from the case where $G$ is not subcubic. If $G$ has a vertex $v$ of degree at least $4$ and none of the components of $G-N[v]$ are $\mathcal{E}_4$-graphs, then we can apply Lemma~\ref{LemIsol2} and the induction hypothesis to obtain a $C_4$-isolating set $D$ of $G-N[v]$ of size at most $n - |N[v]| \leq n-5$, and hence, by Lemma~\ref{LemIsol1} (with $\mathcal{F} = \{C_4\}$, $X = \{v\}$ and $Y = N[v]$), $\iota(G,C_4) \leq n/5$, as required. Clearly, for subcubic graphs, this does not work without further considerations. In this subsection, we prove Theorem~\ref{ThmC4isol} for the case where $G$ is subcubic. 

Each graph in $\cE_4$ has $4$ or $9$ vertices. It is convenient to first settle the case where $G$ has at most $9$ vertices so that we can then focus on the main proof mechanism. 

\begin{lemma}\label{smallsubcubic}
    If $n \leq 9$ and $G$ is a connected subcubic $n$-vertex graph that is not a $\cG_4$-graph, then $\iota(G,C_4) \leq 1$.
\end{lemma}

\noindent
\textbf{Proof.} The result is trivial if $n \leq 4$. Suppose $5 \leq n \leq9$ and $\iota(G,C_4) \geq 1$. Then, $G$ is neither a path nor a cycle, and hence, since $G$ is connected and subcubic, $\Delta(G) = 3$. Let $v \in V(G)$ with $d(v) = 3$. Let $x_1$, $x_2$ and $x_3$ be the neighbours of $v$. Let $G' = G-N[v]$. If $G'$ is $C_4$-free, then $\iota(G,C_4) = 1$. Suppose that $G'$ contains a $4$-cycle $H_1 = (\{y_1, y_2, y_3, y_4\}, \{y_1y_2, y_2y_3, y_3y_4, y_4y_1\})$. Then, $8 \leq n \leq 9$. 

    Suppose $n = 8$. Then, $|V(G')| = 4$. Since $\Delta(G) = 3$ and $G$ is connected, $G'$ is not a $4$-clique, so $G'$ is a $\{C_4,C_4'\}$-graph and $x_i y_j \in E(G)$ for some $i \in [3]$ and $j \in [4]$. We may assume that $i = j = 1$. Since $\Delta(G)=3$, $N[y_1] = \{x_1, y_1, y_2, y_4\}$, so $V(G - N[y_1]) = \{v, x_2, x_3, y_3\}$. Since $y_2, y_4 \in N[y_3]$, $d_{G-N[y_1]}(y_3) \leq \Delta(G) - 2 = 1$. Thus, $G-N[y_1]$ is $C_4$-free, and hence $\iota(G,C_4) = 1$.

    Suppose $n=9$. Let $Y = \{y_1, y_2, y_3, y_4\}$. We have $V(G') = Y \cup \{w\}$ for some $w \in V(G) \setminus Y$. Suppose $x_iy_j \notin E(G)$ for every $i \in [3]$ and $j \in [4]$. Since $G$ is connected, $w x_i, w y_j \in E(G)$ for some $i \in [3]$ and $j \in [4]$. Thus, $G-N[w]$ is $C_4$-free, and hence $\iota(G,C_4)=1$. Now suppose $x_iy_j \in E(G)$ for some $i \in [3]$ and $j \in [4]$. We may assume that $i = j = 1$. Thus, $N[y_1] = \{x_1, y_1, y_2, y_4\}$. 
    
    Let $F_1 = G-N[y_1]$. If $F_1$ is $C_4$-free, then $\iota(G,C_4)=1$. Suppose that $F_1$ contains a $4$-cycle $H_2$. Since $y_2, y_4 \in N[y_3]$, $d_{F_1}(y_3) \leq \Delta(G) - 2 = 1$. Thus, $V(H_2) = \{v, x_2, x_3, w\}$, and hence $wx_2, wx_3 \in E(G)$. Let $F_2 = G-N[x_1]$. If $F_2$ is $C_4$-free, then $\iota(G,C_4)=1$. Suppose that $F_2$ contains a $4$-cycle $H_3$. Let $X = \{x_2, x_3, w\}$ and $Y' = \{y_2, y_3, y_4\}$. We have $V(H_3) \subseteq V(F_2) \subseteq X \cup Y'$. Thus, $xy \in E(H_3)$ for some $x \in X$ and $y \in Y'$. Since $d_{H_2}(x) \geq 2$, we have $|N(x) \cap Y'| \leq 1$, so $xx' \in E(H_3)$ for some $x' \in X \setminus \{x\}$. Since $d_{H_1}(y) \geq 2$, we have $|N(y) \cap X| \leq 1$, so $yy' \in E(H_3)$ for some $y' \in Y' \setminus \{y\}$. Thus, $E(H_3) = \{xy, xx', yy', x'y'\}$. Suppose $w \notin \{x, x'\}$. Then, $\{x, x'\} = \{x_2, x_3\}$. This yields $N(x) = \{v, x', w, y\}$, which contradicts $\Delta(G) = 3$. Therefore, $w \in \{x, x'\}$. We may assume that $w = x'$. Let $y_x, y_w \in \{y, y'\}$ such that $xy_x, wy_w \in E(H_3)$ (so $y_x y_w = yy'$), and let $y''$ be the unique member of $Y' \setminus \{y, y'\}$. Let $F_3 = G - N[w]$. Then, $V(F_3) = \{v, x_1, y_1, y_x, y''\}$. Since $N_{F_3}(v) = \{x_1\}$ and $N_{F_3}(y_x) = N(y_x) \setminus \{x, y_w\}$ (so $d_{F_3}(y_x) = 1$), $|\{z \in V(F_3) \colon d_{F_3}(z) \geq 2\}| \leq 3$. Thus, $F_3$ contains no $4$-cycle, and hence $\iota(G,C_4) = 1$.\qed

\begin{theorem}\label{ThmC4isolsubcubic}
    If $G$ is a connected subcubic $n$-vertex graph that is not a $\cG_4$-graph, then
    \[\iota(G,C_4) \leq \floor*{\frac{n}{5}}.\]
\end{theorem}

\noindent
\textbf{Proof.} We use induction on $n$. If $n\leq9$, then the result is given by Lemma~\ref{smallsubcubic}. Suppose $n \geq 10$. Since $\iota(G,C_4)$ is an integer, it suffices to show that $\iota(G,C_4)\leq n/5$.  The result is trivial if $\iota(G,C_4)=0$, so suppose $\iota(G,C_4) \geq 1$. Then, $G$ contains a $4$-cycle $A = (\{x_0, x_1, x_2, w\}, \{x_0x_1, x_0x_2, x_1w, x_2w\})$. Since $n > 4$ and $G$ is connected, $vx_3 \in E(G)$ for some $v \in V(A)$ and $x_3 \in V(G) \setminus V(A)$. We may assume that $v = x_0$. Since $G$ is subcubic, $N(v)= \{x_1, x_2, x_3\}$ and $\Delta(G) = 3$.  

    Let $G' = G-N[v]$, and let $\cH$ be the set ${\rm C}(G')$ (of components of $G'$). Let $H_w$ be the member of $\cH$ such that $w \in V(H_w)$. Let $\cH' = \{H \in \cH \colon \iota(H,C_4) > |V(H)|/5\}$, and let $\cH^* = \cH \setminus \cH'$. Since $|N(w) \cap N[v]| \geq 2$, we have $d_{H_w}(w) \leq 1$, so $H_w \in \cH^*$. Let $H_w' = H_w - w$.

    Suppose $\cH' = \emptyset$. If $H_w'$ is not a $\cG_4$-graph, then by Lemma~\ref{LemIsol1}, Lemma~\ref{LemIsol2}, Corollary~\ref{LemIsolG'cor} and the induction hypothesis, 
\begin{align*}
\iota(G,C_4) &\leq |\{v\}| + \iota(G',C_4) = 1 + \iota(H_w,C_4) + \sum_{H\in\cH\setminus \{H_w\}}\iota(H,C_4) \\
&\leq \frac{|N[v]|+1}{5} + \iota(H_w',C_4) + \sum_{H\in\cH\setminus \{H_w\}} {\frac{|V(H)|}{5}} \\
&\leq \frac{|N[v]|+1}{5}  + \frac{|V(H_w)|-1}{5} +\sum_{H\in\cH\setminus \{H_w\}} {\frac{|V(H)|}{5}} = \frac{n}{5}.
\end{align*}
Suppose that $H_w'$ is a $\cG_4$-graph. Then, $H_w'$ contains a $4$-cycle $A'$, and $V(H_w) = \{w\} \cup V(A')$. Let $w_1w_2$, $w_2w_3$, $w_3w_4$ and $w_4w_1$ be the edges of $A'$. Since $H_w$ is connected, we may assume that $ww_1 \in E(H)$. Thus, $N(w_1) = \{w, w_2, w_4\}$. Let $F_1 = G - N[w_1]$, and let $F_1' = F_1 - w_3$. Thus, $F_1' = G - V(H_w)$. Since $n \geq 10$, $|V(F_1')| \geq 5$. Clearly, $F_1'$ is connected (as $H_w \in \mathcal{H}$, $N[v] \subseteq V(F_1')$, and, since $G$ is connected, each member of $\mathcal{H}$ has a vertex that is adjacent to a member of $N(v)$), so $F_1'$ is not a $\cG_4$-graph. Since $w_2, w_4 \in N(w_3)$, $d_{F_1}(w_3) \leq 1$. By Corollary~\ref{LemIsolG'cor}, $\iota(F_1,C_4) = \iota(F_1',C_4)$.  By Lemma~\ref{LemIsol1} and the induction hypothesis,
\[\iota(G,C_4) \leq 1 + \iota(F_1,C_4) = \frac{|N[w_1]| + 1}{5} + \iota(F_1',C_4) \leq \frac{|N[w_1]| + 1}{5} + \frac{|V(F_1)| - 1}{5} = \frac{n}{5}. \]

    Now suppose $\cH' \neq \emptyset$. By the induction hypothesis, each member of $\cH'$ is a $\cG_4$-graph. For each $H \in \cH^*$, let $D_H$ be a $C_4$-isolating set of $H$ of size $\iota(H,C_4)$. By the definition of $\mathcal{H}^*$, $|D_H| \leq |V(H)|/5$ for each $H \in \mathcal{H}^*$.
    
    For any $H \in \cH$ and $x \in N(v)$ such that $xy_{x,H} \in E(G)$ for some $y_{x,H} \in V(H)$, we say that $H$ is \emph{linked to $x$} and that $x$ is \emph{linked to $H$}. Since $G$ is connected, each member of $\cH$ is linked to at least one member of $N(v)$. For each $x \in N(v)$, let $\cH_{x} = \{H \in \cH \colon H$ is linked to $ x\}$, $\cH'_{x} = \{H \in \cH' \colon H $ is linked to $x\}$ and $\cH_{x}^* = \{H \in \cH^* \colon H$ is linked to $x$ only$\}$. 
    
    We have 
\begin{equation} n = |N[v]| + 4|\mathcal{H}'| + \sum_{H \in \mathcal{H}^*} |V(H)| \geq |N[v]| + 4|\mathcal{H}_x'| + 4|\mathcal{H}' \setminus \mathcal{H}_x'| + \sum_{H \in \mathcal{H}^*} 5|D_H|. \label{ineq for n}
\end{equation}

\noindent
\textbf{Case 1}: \textit{$|\cH_x'| \geq 2$ for some $x\in N(v)$.} 
For each $H \in \cH'\setminus\cH_x'$, let $x_H \in N(v)$ such that $H$ is linked to $x_H$. Let $X = \{x_H \colon H \in \cH'\setminus \cH'_x\}$. Note that $x \notin X$. Let $D = \{v,x\} \cup X \cup \bigcup_{H\in \cH^*} D_H$. We have $y_{x,H} \in N[x]$ for each $H \in \cH_x'$, and $y_{x_H,H} \in N[x_H]$ for each $H \in \cH' \setminus \cH'_{x}$, so $D$ is a $C_4$-isolating set of $G$. By (\ref{ineq for n}) and $|\cH_x'| \geq 2$, 
\[n \geq |\{v, x\} \cup X| + 8 + 4|X| + \sum_{H \in \mathcal{H}^*} 5|D_H| = 10 + 5|X| + \sum_{H \in \mathcal{H}^*} 5|D_H| \geq 5|D|.\]
We have $\iota(G, C_4) \leq |D| \leq n/5$.\medskip

\noindent
\textbf{Case 2}: \textit{$|\cH'_{x}| \leq 1$ for each $x \in N(v)$.} For each $H \in \cH'$, let $x_H \in N(v)$ such that $H$ is linked to $x_H$.\medskip
    
\noindent
\textbf{Subcase 2.1}: \textit{For some $H \in \cH'$, $H$ is linked to $x_H$ only.} 

Suppose $x_H \in \{x_1, x_2\}$. We may assume that $x_H = x_1$. Let $G^* = G - (\{x_1\} \cup V(H))$. Recall that $x_1w, x_2w \in E(G)$. Since $H_w \notin \cH'$, we have $H_w \neq H$, so $w \neq y_{x_1,H}$. Since $\Delta(G) = 3$, $N(x_1) = \{v, w, y_{x_1,H}\}$. Thus, $G^*$ is connected. Now $H$ is a $4$-vertex graph containing a $4$-cycle $(\{u_1, u_2, u_3, u_4\}, \{u_1u_2, u_2u_3, u_3u_4, u_4u_1\})$ with $u_1 = y_{x_1,H}$. Since $n \geq 10$, $G^*$ is not a $\cG_4$-graph. We have $N(y_{x_1,H}) = \{x_1, u_2, u_4\}$. Let $F_1 = G - N[y_{x_1,H}]$. Then, $G^* = F_1 - u_3$. Since $H$ is linked to $x_1$ only, we have $d_{F_1}(u_3) = 0$, so $\iota(F_1,C_4) = \iota(G^*,C_4)$. By Lemma~\ref{LemIsol1} and the induction hypothesis, $\iota(G,C_4) \leq 1 + \iota(F_1,C_4) = 1 + \iota(G^*,C_4) \leq 1 + (n-5)/5 = n/5$.

Now suppose $x_H = x_3$. Let $G^* = G - (\{x_3\} \cup V(H))$. Then, $G^*$ has a component $G_v^*$ with $(N[v] \setminus \{x_3\}) \cup V(H_w) \subseteq V(G_v^*)$, and any other component of $G^*$ is a member of $\cH_{x_3}^*$ (as we are in Case~2). Let $D^*$ be a $C_4$-isolating set of $G_v^*$ of size $\iota(G_v^*,C_4)$. Let $D' = D^* \cup \{x_3\} \cup \bigcup_{I \in \cH_{x_3}^*} D_{I}$. Then, $D'$ is a $C_4$-isolating set of $G$. Thus, $\iota(G,C_4) \leq |D^*| + 1 + \sum_{I \in \cH_{x_3}^*} |D_{I}|$, and hence
\[\iota(G,C_4) \leq |D^*| + \frac{|\{x_3\} \cup V(H)|}{5} + \sum_{I \in \cH_{x_3}^*} \frac{|V(I)|}{5} = |D^*| + \frac{n-|V(G_v^*)|}{5}.\]
Thus, $\iota(G,C_4) \leq n/5$ if $|D^*| \leq |V(G_v^*)|/5$. Suppose $|D^*| > |V(G_v^*)|/5$. By the induction hypothesis, $G_v^*$ is a $\cG_4$-graph. Since $v \in N[x_3]$, $D' \setminus D^*$ is a $C_4$-isolating set of $G$, so $\iota(G,C_4) \leq (n-|V(G_v^*)|)/5 < n/5$.\medskip

\noindent
\textbf{Subcase 2.2}: \textit{For each $H \in \cH'$, $H$ is linked to at least $2$ distinct neighbours of $v$.} Let $H \in \cH'$. For some $i, i' \in [3]$ with $i < i'$, $H$ is linked to $x_i$ and $x_{i'}$. Since we are in Case~2, for each $j \in \{i, i'\}$, $H$ is the only member of $\cH'$ that is linked to $x_j$. Thus, any member of $\cH' \setminus \{H\}$ can only be linked to the unique member of $N(v) \setminus \{x_i, x_{i'}\}$. Since we are in Subcase~2.2, $\cH' \setminus \{H\} = \emptyset$, so $\cH' = \{H\}$.  We have $i \in \{1, 2\}$. Let $G^* = G - (\{x_i\} \cup V(H))$. Then, $N(x_i) = \{v, w, y_{x_i,H}\}$ and $G^*$ is connected. Now $H$ is a $4$-vertex graph containing a $4$-cycle $(\{u_1, u_2, u_3, u_4\}, \{u_1u_2, u_2u_3, u_3u_4, u_4u_1\})$ with $u_1 = y_{x_i,H}$. Since $n \geq 10$, $G^*$ is not a $\cG_4$-graph. We have $N(y_{x_i,H}) = \{x_i, u_2, u_4\}$. Let $F_1 = G - N[y_{x_i,H}]$. Then, $G^* = F_1 - u_3$. Since $u_2, u_4 \in N(u_3)$, we have $d_{F_1}(u_3) \leq 1$, so $\iota(F_1,C_4) = \iota(G^*,C_4)$ by Corollary~\ref{LemIsolG'cor}. By Lemma~\ref{LemIsol1} and the induction hypothesis, $\iota(G,C_4) \leq 1 + \iota(F_1,C_4) = 1 + \iota(G^*,C_4) \leq 1 + (n-5)/5 = n/5$.\qed

\subsection{The general case}

We now use Theorem~\ref{ThmC4isolsubcubic} to prove Theorem~\ref{ThmC4isol} in its entirety. By Theorem~\ref{ThmC4isolsubcubic}, we now need to consider graphs whose maximum degree is at least $4$. As in Section~\ref{subcubic section}, it is convenient to first settle the case where $n \leq 9$. 

\begin{lemma}\label{ThmC4isolsmall}
    If $n \leq 9$ and $G$ is a connected $n$-vertex graph that is not an $\cE_4$-graph, then $\iota(G,C_4) \leq 1$.
\end{lemma}

\noindent
\textbf{Proof.} The result is trivial if $n\leq 4$. If $\Delta(G) \leq 3$, then the result is given by Theorem \ref{ThmC4isolsubcubic}. Suppose $5 \leq n \leq 9$ and $\Delta(G) \geq 4$. Let $v \in V(G)$ such that $d(v) = \Delta(G)$. If either $n \leq 8$ or $n = 9$ and $\Delta(G) \geq 5$, then $G - N[v]$ has at most $3$ vertices, and hence $\{v\}$ is a $C_4$-isolating set of $G$. It remains to consider $n = 9$ and $\Delta(G) = 4$. The appendix contains a code, written in Mathematica by the first author, that performs an exhaustive search for graphs which violate the bound in Theorem~\ref{ThmC4isol}. It determines that the $\cG_9$-graphs are the connected subquartic $9$-vertex graphs which violate the bound. \qed 
\medskip

The following are properties of $\cG_9$-graphs that we use in the proof of Theorem~\ref{ThmC4isol}.

\begin{observation}\label{ObsG9} Let $G \in \mathcal{G}_9$. \\
(a) $\iota(G,C_4) = 2$.\\
(b) $\delta(G) \geq 3$ and $\Delta(G) = 4$. Moreover, $G$ is $4$-regular if it is not one of $G_5$ and $G_6$.\\
(c) If $G$ is $G_5$ or $G_6$, then no two vertices of $G$ of degree $3$ are adjacent.\\
(d) For each $v \in V(G)$, $G - (\{v\} \cup N[v'])$ is a $\{P_3,K_3\}$-graph for some $v' \in V(G)\setminus \{v\}$.
\end{observation}

\noindent
\textbf{Proof of Theorem \ref{ThmC4isol}.} We first show that the bound is attained if $G$ is the graph $B_{n,C_4}$ in Construction~\ref{Bconstruction}. Let $B = B_{n,C_4}$. If $n\leq 4$, then $\iota(B,C_4) = 0$. Suppose $n\geq 5$. Since $[a_{n,4}]$ is a $C_4$-isolating set of $B$, $\iota(B,C_4) \leq a_{n,4}$. If $D$ is a $C_4$-isolating set of $B$, then $D \cap (V(F_i) \cup \{i\}) \neq \emptyset$ for each $i \in [a_{n,4}]$, so $\iota(B, C_4) \geq a_{n,4}$. Together with $\iota(B,C_4) \leq a_{n,4} =  \big\lfloor \frac{n}{5} \big\rfloor$, this gives us $\iota(B,C_4) =  \big\lfloor \frac{n}{5} \big\rfloor$.

    We now prove that the bound holds, using induction on $n$. It suffices to prove that $\iota(G,C_4) \leq n/5$ as $\iota(G,C_4)$ is an integer. If $n \leq 9$, then the bound holds by Lemma~\ref{ThmC4isolsmall}. Suppose $n \geq 10$. If $\Delta(G) \leq 3$, then the result is given by Theorem~\ref{ThmC4isolsubcubic}. Suppose $\Delta(G) \geq 4$. Let $v \in V(G)$ such that $d(v) = \Delta(G)$. Then, $|N[v]| \geq 5$. Let $G' = G-N[v]$. Then, $|V(G')| \leq n-5$. If $V(G') = \emptyset$, then $\iota(G, C_4) \leq 1 \leq n/5$. Suppose $V(G') \neq \emptyset$. Let $\cH = {\rm C}(G')$. Let $\cH' = \{H \in \cH \colon \iota(H, C_4) > |V(H)|/5\}$. 
    
    Suppose $\cH' = \emptyset$. By Lemmas~\ref{LemIsol1} and \ref{LemIsol2},
    \begin{equation*}
        \iota(G,C_4) \leq 1 + \iota(G',C_4) = 1 + \sum_{H\in\cH}\iota(H,C_4) \leq \frac{|N[v]|}{5} + \sum_{H \in \cH}\frac{|V(H)|}{5} = \frac{n}{5}.
    \end{equation*}

    Now suppose $\cH' \neq \emptyset$. By the induction hypothesis, each member of $\cH'$ is an $\cE_4$-graph. For each $H \in \cH \setminus \cH'$, let $D_H$ be a $C_4$-isolating set of $H$ of size $\iota(H, C_4)$. By the definition of $\mathcal{H}'$, $|D_H| \leq |V(H)|/5$ for each $H \in \mathcal{H} \setminus \mathcal{H}'$.

    For each $x \in N(v)$, we define $\mathcal{H}_x$, $\mathcal{H}_x'$ and $\mathcal{H}_x^*$ as in the proof of Theorem~\ref{ThmC4isolsubcubic}. Since $G$ is connected, each member of $\cH$ is linked to at least one member of $N(v)$.  

    For any $x \in N(v)$ and $H \in \cH_x'$, let $H_x' = H - y_{x,H}$. Note that $H_x'$ is connected and subquartic. If $H$ is a $\cG_4$-graph, then we take $D_{x,H} = \emptyset$. If $H$ is a $\cG_9$-graph, then by Observation \ref{ObsG9}, there exists some $y' \in V(H)$ such that $H_x' - N[y']$ is a $\{P_3,K_3\}$-graph, and we take $D_{x,H} = \{y'\}$. Thus, $D_{x,H}$ is a $C_4$-isolating set of $H_x'$, and $|D_{x,H}| \leq (|V(H_x')| - 3)/5$.\medskip

\noindent
\textbf{Case 1}: \textit{$|\cH'_x| \geq 2$ for some $x \in N(v)$.} For each $H \in \cH' \setminus \cH'_x$, let $x_H \in N(v)$ such that $H$ is linked to $x_H$. Let $X = \{x_H \colon H \in \cH' \setminus \cH'_x\}$. Note that $x \notin X$. 
    Let 
    \[D = \{v,x\} \cup X \cup \left( \bigcup_{H \in \cH_x'} D_{x,H} \right) \cup \left( \bigcup_{H \in \cH' \setminus \cH_x'} D_{x_H,H} \right) \cup \left( \bigcup_{H \in \cH \setminus \cH'} D_{H} \right).\]
We have $y_{x,H} \in N[x]$ for each $H \in \cH_x'$, and $y_{x_H,H} \in N[x_H]$ for each $H \in \cH' \setminus \cH'_{x}$, so $D$ is a $C_4$-isolating set of $G$. Let 
    \[Y = \{v, x\} \cup X \cup \{y_{x,H} \colon H \in \cH_x'\} \cup \{y_{x_H,H} \colon H \in \cH' \setminus \cH'_x\}.\]
Since $|\cH'_x|\geq 2$ and $|\cH'\setminus\cH'_x| \geq |X|$,  $|Y| \geq 2 + |X| + 2 + |X| \geq 10 + 5|X| - 3|\cH'_x| - 3|\cH'\setminus\cH'_x|$. We have
\begin{align*}
\iota(G, C_4) &\leq |D| = 2 + |X| + \sum_{H \in \cH_x'} |D_{x,H}| + \sum_{H \in \cH' \setminus \cH_x'} |D_{x_H,H}| + \sum_{H \in \cH \setminus \cH'} |D_{H}| \\
&\leq \frac{10 + 5|X|}{5} + \sum_{H \in \cH_x'} \frac{|V(H_x')| - 3}{5} + \sum_{H \in \cH' \setminus \cH_{x}'} \frac{|V(H_{x_H}')| - 3}{5} + \sum_{H \in \cH \setminus \cH'} \frac{|V(H)|}{5} \\
&\leq \frac{|Y|}{5} + \sum_{H \in \cH_x'} \frac{|V(H_x')|}{5} + \sum_{H \in \cH' \setminus \cH_{x}'} \frac{|V(H_{x_H}')|}{5} + \sum_{H \in \cH \setminus \cH'} \frac{|V(H)|}{5} \leq \frac{n}{5}.
\end{align*}

\noindent
\textbf{Case 2}:
\begin{equation} |\cH'_x| \leq 1 \mbox{ \emph{for each} } x \in N(v). \label{EqHxleq1} 
\end{equation} 
For each $H \in \cH'$, let $x_H \in N(v)$ such that $H$ is linked to $x_H$. Let $X = \{x_H \colon H \in \cH'\}$. Since no two members of $\cH'$ are linked to the same neighbour of $v$, $|X| = |\cH'|$. Let $W = N(v)\setminus X$, $W' = N[v]\setminus X$ and $Y_X=X\cup\{y_{x_H,H} \colon H\in\cH'\}$. Note that $|Y_X| = |X| + |\cH'| = 2|X|$.\medskip

\noindent
\textbf{Subcase 2.1}: $|W| \geq 4$. Let $D = \{v\} \cup X \cup \left( \bigcup_{H \in \cH'} D_{x_H,H} \right) \cup ( \bigcup_{H \in \cH \setminus \cH'} D_{H} )$. Then, $D$ is a $C_4$-isolating set of $G$. Similarly to Case~1, we obtain
\begin{align*} \iota(G, C_4) &\leq |D| = 1 + |X| +  \sum_{H \in \cH'} |D_{x_H,H}| + \sum_{H \in \cH \setminus \cH'} |D_{H}|\\
&\leq \frac{|\{v\}\cup W| + 5|X|}{5} + \sum_{H\in\cH'}\frac{|V(H_{x_H}')|-3}{5} + \sum_{H\in\cH\setminus\cH'} \frac{|V(H)|}{5} \\
&= \frac{|W'|+|Y_X|}{5} +\sum_{H\in\cH'}\frac{|V(H_{x_H}')|}{5} + \sum_{H\in\cH\setminus\cH'} \frac{|V(H)|}{5} = \frac{n}{5}.
\end{align*}

\noindent
\textbf{Subcase 2.2}: $|W| \leq 3$.\medskip
    
\noindent
\textbf{Subcase 2.2.1}: \textit{For some $H \in\cH'$, $H$ is linked to $x_H$ only.} Let $G^* = G - (\{x_H\} \cup V(H))$. Then, $G^*$ has a component $G_v^*$ such that $N[v] \setminus \{x_H\} \subseteq V(G_v^*)$, and the other components of $G^*$ (if any) are the members of $\cH_{x_H}^*$. Let $D^*$ be a $C_4$-isolating set of $G_v^*$ of size $\iota(G_v^*, C_4)$. Let $D = D^* \cup \{x_H\} \cup D_{x_H,H} \cup \bigcup_{I \in \cH_{x_H}^*} D_I$. Since $D$ is a $C_4$-isolating set of $G$,
\begin{equation} \iota(G, C_4) \leq |D^*| + 1 + |D_{x_H,H}| + \sum_{I \in \cH_{x_H}^*} |D_I| \leq \iota(G_v^*, C_4) + \frac{|\{x_H\} \cup V(H)|}{5} + \sum_{I \in \cH_{x_H}^*} \frac{|V(I)|}{5}.  \nonumber
\end{equation}
This yields $\iota(G, C_4) \leq n/5$ if $\iota(G_v^*, C_4) \leq |V(G_v^*)|/5$. Suppose $\iota(G_v^*, C_4) > |V(G_v^*)|/5$. By the induction hypothesis, $G_v^*$ is an $\cE_4$-graph. If $|V(G_v^*)| = 4$, then, since $v \in N[x_H]$, $D \setminus D^*$ is a $C_4$-isolating set of $G$, so $\iota(G,C_4) \leq |\{x_H\} \cup V(H)|/5 + \sum_{I \in \cH_{x_H}^*} |V(I)|/5 < n/5$. Suppose $|V(G_v^*)| = 9$. By Observation~\ref{ObsG9}, $G_v^*$ has a vertex $w$ such that $G_v^* - (\{v\} \cup N[w])$ is a $\{P_3,K_3\}$-graph. Since $v \in N[x_H]$, $(D \setminus D^*) \cup \{w\}$ is a $C_4$-isolating set of $G$, so
\begin{equation*} \iota(G,C_4) \leq 1+\frac{|\{x_H\} \cup V(H)|}{5} + \sum_{I \in \cH_{x_H}^*} \frac{|V(I)|}{5} < \frac{|V(G_v^*)|}{5} + \frac{n-|V(G_v^*)|}{5} = \frac{n}{5}.
\end{equation*}

\noindent
\textbf{Subcase 2.2.2}: \textit{For each $H \in \cH'$, $H$ is linked to some $x_H' \in N(v) \setminus \{x_H\}$.} Let $X' = \{x_H' \colon H \in \cH'\}$. By (\ref{EqHxleq1}), for every $H, I \in \mathcal{H}'$ with $H \neq I$, we have $x_H' \neq x_I'$ and $x_H \neq x_I'$. Thus, $|X'| = |\cH'|$ and $X' \subseteq W$. This yields $|\cH'| \leq 3$ since $|W| \leq 3$. 

Let $h = |\mathcal{H}'|$. We have $1 \leq h \leq 3$. Let $H_1, \dots, H_h$ be the members of $\mathcal{H}'$. For each $i \in [h]$, let $x_i = x_{H_i}$ and $x_i' = x_{H_i}'$. Let $Y_1 = N[y_{x_1, H_1}] \cap (\{x_1\} \cup V(H_1))$. Let $G^* = G - Y_1$, $Y_1' = V(H_1) \setminus Y_1$ and $H_1^* = H_1[Y_1']$ ($= H_1 - N_{H_1}[y_{x_1, H_1}]$).\medskip

\noindent
\textbf{Subcase 2.2.2.1}: $h = 3$. Then, $|W| = 3$, $X' = W$, $N(v) = X \cup W = \{x_1, x_2, x_3, x_1', x_2', x_3'\}$, and $G^*$ has a component $G_v^*$ such that $N[v] \setminus \{x_1\}, V(H_2), V(H_3) \subseteq V(G_v^*)$. Since $d(v) = 6$ and $|V(H_i)| \in \{4, 9\}$ for each $i \in [3]$, we have $|V(G_v^*)| \notin \{4, 9\}$, so $G_v^*$ is not an $\cE_4$-graph. By (\ref{EqHxleq1}), $N[V(H_1)] \cap \{x_2, x_3, x_2', x_3'\} = \emptyset$, so
\begin{equation} N_{G^*}[V(H_1^*)] \cap N(v) \subseteq \{x_1'\}. \label{H1links1}
\end{equation}

Suppose $|V(H_1)| = 4$. If $Y_1 = \{x_1\} \cup V(H_1)$, then ${\rm C}(G^*) = \{G_v^*\} \cup \cH_{x_1}^*$, and by Lemma~\ref{LemIsol1}, Lemma~\ref{LemIsol2} and the induction hypothesis, 
\[\iota(G,C_4) \leq 1 + \iota(G_v^*,C_4) + \sum_{H \in \cH_{x_1}^*} \iota(H,C_4) \leq \frac{|Y_1|}{5} + \frac{|V(G_v^*)|}{5} + \sum_{H \in \cH_{x_1}^*} \frac{|V(H)|}{5} = \frac{n}{5}.\]
Suppose $Y_1 \neq \{x_1\} \cup V(H_1)$. Then, $Y_1 = (\{x_1\} \cup V(H_1)) \setminus \{y^*\}$ with $\{y^*\} = V(H_1^*)$. Let $G^- = G^* - y^*$, let $G_v^- = G_v^* - y^*$ if $y^* \in V(G_v^*)$, and let $G_v^- = G_v^*$ if $y^* \notin V(G_v^*)$. Then, ${\rm C}(G^-) = \{G_v^-\} \cup \cH_{x_1}^*$. Since $|V(G_v^-)| \notin \{4, 9\}$, $G_v^-$ is not an $\cE_4$-graph. By (\ref{H1links1}) and Corollary~\ref{LemIsolG'cor}, $\iota(G^*,C_4) = \iota(G^-,C_4)$. By Lemma~\ref{LemIsol1}, $\iota(G,C_4) \leq 1 + \iota(G^*,C_4) = 1 + \iota(G^-,C_4)$. Therefore, by Lemma~\ref{LemIsol2} and the induction hypothesis, 
\begin{equation}\label{EqD-} \iota(G,C_4) \leq 1 + \iota(G_v^-,C_4) + \sum_{H \in \cH_{x_1}^*} \iota(H,C_4) \leq \frac{|Y_1 \cup \{y^*\}|}{5} + \frac{|V(G_v^-)|}{5} + \sum_{H \in \cH_{x_1}^*} \frac{|V(H)|}{5} = \frac{n}{5}. \nonumber
\end{equation}

Now suppose $|V(H_1)| = 9$. By Observation~\ref{ObsG9}, we have $\iota(H_1,C_4) = 2$, $\delta(H_1) \geq 3$ and $\Delta(H_1) = 4$. Also, $G_5$ and $G_6$ are the only members of $\cG_9$ whose minimum degree is $3$, and we see that they have exactly $2$ and $4$ vertices of degree $3$, respectively, no two of which are adjacent. Since $\iota(H_1,C_4) = 2$, $H_1-Y_1$ contains a $4$-cycle. It follows that $H_1-Y_1$ has at least $2$ vertices $z$ and $z'$ such that $d_{H_1}(z) = d_{H_1}(z') = 4$. Let $Z = N_{H_1}[z]$. Since $\iota(H_1,C_4) = 2$, $H_1-Z$ contains a $4$-cycle $A_Z$. We have $|V(H_1-Z)| = 4$, so $V(H_1-Z) = V(A_Z)$. Since $z \notin Y_1$, we have $z \notin N[y_{x_1,H_1}]$, so $y_{x_1,H_1} \in V(A_Z)$, and hence $G-Z$ is connected. Since $|V(G-Z)| > 9$, $\iota(G-Z,C_4) \leq (n-5)/5$ by the induction hypothesis. By Lemma~\ref{LemIsol1}, $\iota(G,C_4) \leq 1 + \iota(G-Z,C_4) \leq n/5$.\medskip

\noindent
\textbf{Subcase 2.2.2.2}: $h = 2$. Then, $G^*$ has a component $G_v^*$ such that $N[v] \setminus \{x_1\}, V(H_2) \subseteq V(G_v^*)$. We have $4 \leq \Delta(G) = d(v) = |X| + |W| = h + |W| = 2 + |W|$, so $|W| \geq 2$. Recall that $|W| \leq 3$.   

Suppose $|W| = 3$. Then, $N(v) = \{x_1,x_1',x_2,x_2',w\}$ for some $w \in V(G) \setminus \{x_1,x_1',x_2,x_2'\}$. By (\ref{EqHxleq1}), each vertex in $N(v)$ is linked to at most one of $H_1$ and $H_2$. We may assume that $w$ is not linked to $H_1$. Thus, $N[V(H_1)] \setminus V(H_1) \subseteq \{x_1, x_1'\}$. If either $|V(H_1)| = 4$ and $G_v^* - (Y_1' \cap V(G_v^*))$ is not a $\cG_9$-graph, or $|V(H_1)| = 9$, then we obtain $\iota(G,C_4) \leq n/5$ as in Subcase~2.2.2.1. Suppose that $|V(H_1)| = 4$ and $G_v^* - (Y_1' \cap V(G_v^*))$ is a $\cG_9$-graph. Then, $|V(H_2)| = 4$. We have $|N(x) \cap V(H_1)| \geq 1$ for each $x \in \{x_1, x_1'\}$.

Suppose $|N(x) \cap V(H_1)| \geq 2$ for some $x \in \{x_1, x_1'\}$. Let $x'$ be the member of $\{x_1, x_1'\} \setminus \{x\}$, and let $z_1, z_2, z_3$ and $z_4$ be the $4$ vertices of $H_1$, where $z_1, z_2 \in N(x)$. Let $Y = \{z_3, z_4\}$, $Z_1 = \{x, z_1, z_2\}$ and $Z_2 = Z_1 \cup \{v\}$. For each $i \in [2]$, let $G_{Z_i} = G - Z_i$ and $F_i = G_{Z_i} - Y$. Then, $F_1 = G - (\{x\} \cup V(H_1))$, $F_2 = F_1 - v$, $F_1$ has a component $F_v^*$ such that $N[v] \setminus \{x\}, V(H_2) \subseteq V(F_v^*)$, and ${\rm C}(F_1) = \{F_v^*\} \cup \cH_x^*$. By Lemma~\ref{LemIsolG'}, $\iota(G_{Z_i},C_4) = \iota(F_i,C_4)$ as $x' \in V(G_{Z_i})$, $N_{G_{Z_i}}[Y] \cap V(F_i) \subseteq \{x'\}$ and $G_{Z_i}[\{x'\} \cup Y]$ contains no $4$-cycle. We have $|V(F_v^*)| \geq d(v) + |V(H_2)| > 4$. Suppose that $F_v^*$ is not a $\cG_9$-graph. By Lemma~\ref{LemIsol2} and the induction hypothesis, $\iota(F_1,C_4) \leq |V(F_v^*)|/5 + \sum_{H \in \cH_x^*} |V(H)|/5 = |V(F_1)|/5$. Since $Z_1 \subset N[x]$, $\iota(G,C_4) \leq 1 + \iota(G_{Z_1},C_4)$ by Lemma~\ref{LemIsol1}. Since $\iota(G_{Z_1},C_4) = \iota(F_1,C_4)$, we obtain $\iota(G,C_4) \leq 1 + \iota(F_1,C_4) \leq 1 + |V(F_1)|/5 = n/5$. Now suppose that $F_v^*$ is a $\cG_9$-graph. Note that $F_2$ has a component $F_2'$ such that $\{x_2, x_2'\} \cup V(H_2) \subseteq V(F_2') \subseteq V(F_v^*-v)$, and ${\rm C}(F_2) = {\rm C}(F_v^*-v) \cup \cH_x^*$. Thus, $6 \leq |V(F_2')| \leq |V(F_v^*-v)| = 8$. By the induction hypothesis, $\iota(F_2',C_4) = 1$ and $\iota(I,C_4) = 0$ for any $I \in {\rm C}(F_v^*-v) \setminus \{F_2'\}$. Thus, by Lemma~\ref{LemIsol2}, $\iota(F_2,C_4) \leq 1 + \sum_{H \in \cH_x^*} |V(H)|/5 < |V(F_2)|/5$.
Since $Z_2 \subseteq N[x]$, $\iota(G,C_4) \leq 1 + \iota(G_{Z_2},C_4)$ by Lemma~\ref{LemIsol1}. Since $\iota(G_{Z_2},C_4) = \iota(F_2,C_4)$, we obtain $\iota(G,C_4) \leq 1 + \iota(F_2,C_4) < 1 + |V(F_2)|/5 < n/5$.

Now suppose $|N(x) \cap V(H_1)| = 1$ for each $x \in \{x_1, x_1'\}$. Let $x = x_1$ and $x' = x_1'$. Let $z_1$ be the member of $N(x) \cap V(H_1)$. Let $Y = V(H_1) \setminus \{z_1\}$, $Z_1 = \{x, z_1\}$ and $Z_2 = Z_1 \cup \{v\}$. For each $i \in [2]$, let $G_{Z_i} = G - Z_i$ and $F_i = G_{Z_i} - Y$. Then, $F_1 = G - (\{x\} \cup V(H_1))$ and $F_2 = F_1 - v$. Since $|\{x'\} \cup Y| = 4$ and $|N(x') \cap Y| \leq |N(x') \cap V(H_1)| = 1$, $G_{Z_i}[\{x'\} \cup Y]$ contains no $4$-cycle. Thus, we can proceed as in the case $|N(x) \cap V(H_1)| \geq 2$ to obtain $\iota(G,C_4) \leq n/5$.

If $|W| = 2$, then by applying the argument for the case $|W| = 3$, we again obtain $\iota(G, C_4) \leq n/5$.\medskip

\noindent
\textbf{Subcase 2.2.2.3}: $h = 1$. We have $4 \leq \Delta(G) = d(v) = |X| + |W| \leq 1 + 3 =4$, so $\Delta(G) = 4 = d(v)$. Let $Y = N[v] \cup V(H_1)$.  Let $y_1 = y_{x_1,H_1}$ and $y_1' = y_{x_1',H_1}$.

Suppose $|V(H_1)|=9$. Since $\Delta(G) = 4$ and $x_1y_1, x_1'y_1' \in E(G)$, we have $d_{H_1}(y_1) = d_{H_1}(y_1') = 3$ (by Observation~\ref{ObsG9}), so $H_1$ is a $\{G_5, G_6\}$-graph. Suppose $H_1 \simeq G_5$. We may assume that $H_1 = G_5$. Thus, $y_1 \in \{2, 3\}$. Since $d(9) \leq \Delta(G) = 4 = d_{H_1}(9)$, $N(9) = N_{H_1}(9) = \{4,5,6,8\}$. Let $G^* = G - N[9]$. Then, $G^*$ is connected and $|V(G^*)| \geq 1 + d(v) + |V(H_1 - N[9])| = 9$. We have $d(1) \leq \Delta(G) = 4 = d_{H_1}(1)$, so $N(1) = N_{H_1}(1) = \{2,3,4,5\}$, and hence $d_{G^*}(1) = 2$. By Observation~\ref{ObsG9}, $G^*$ is not a $\cG_9$-graph. Thus, $G^*$ is not an $\cE_4$-graph. By Lemma~\ref{LemIsol1} and the induction hypothesis, $\iota(G, C_4) \leq 1 + \iota(G^*,C_4) \leq |N[9]|/5 + |V(G^*)|/5 = n/5$. Now suppose $H_1 \simeq G_6$. We may assume that $H_1 = G_6$. We now obtain $y_1 \in \{2, 3, 4, 5\}$. If $y_1 \in \{2, 3\}$, then the same argument applies. If $y_1 \in \{4, 5\}$, then we can apply the same argument by taking $G^* = G - N[6]$ instead.

Now suppose $|V(H_1)| = 4$. Since $H_1$ is a $\cG_4$-graph, $H_1$ contains a $4$-cycle $H_1'$. Thus, $V(H_1) = V(H_1') = \{y_1, \dots, y_4\}$ and $E(H_1') = \{y_1y_2, y_2y_3, y_3y_4, y_4y_1\} \subseteq E(H_1)$ for some $y_2, y_3, y_4 \in V(G) \setminus \{y_1\}$. Since $|N[v] \cup V(H_1)| = 9 < n$, $\cH \setminus\cH' \neq \emptyset$. Since $\Delta(G) = 4$, $|N(y) \cap N(v)| \leq 2$ for each $y \in V(H_1)$.

Suppose $|N(y) \cap N(v)| = 2$ for some $y \in V(H_1)$. We may assume that $y = y_1$. We have $N(y_1) \cap N(v) = \{x_1, x_1^*\}$ for some $x_1^* \in N(v) \setminus \{x_1\}$. Let 
\[\mathcal{I} = \{H \in \cH \colon H \mbox{ is only linked to one or both members of } \{x_1, x_1^*\}\}.\] 
Let $G^* = G - N[y_1]$. Then, $G^*$ has a component $G_v^*$ such that $N[v] \setminus \{x_1, x_1^*\} \subseteq V(G_v^*)$, the members of $\mathcal{I}$ are components of $G^*$, and if $G^*$ has another component $A^*$, then $V(A^*) = \{y_3\}$ (so $\iota(A^*,C_4) = 0$). By Lemma~\ref{LemIsol1}, $\iota(G,C_4) \leq 1 + \iota(G^*, C_4)$. Thus, by Lemma~\ref{LemIsol2},
\begin{equation} \iota(G,C_4) \leq 1 + \iota(G_v^*, C_4) + \sum_{H \in\mathcal{I}} \iota(H, C_4) \leq \frac{|N[y_1]|}{5} + \iota(G_v^*, C_4) + \sum_{H \in \mathcal{I}} \frac{|V(H)|}{5}. \label{h'=1a}
\end{equation}
This gives $\iota(G,C_4) \leq n/5$ if $\iota(G_v^*,C_4) \leq |V(G_v^*)|/5$. Suppose $\iota(G_v^*,C_4) > |V(G_v^*)|/5$. By the induction hypothesis, $G_v^*$ is an $\cE_4$-graph. By Observation \ref{ObsG9}, $G_v^*$ is not a $\cG_9$-graph as $d_{G_v^*}(v) = 2$. Thus, $|V(G_v^*)| = 4$. Let $w$ and $w'$ be the two members of $N_{G_v^*}(v)$. Since $G_v^*$ contains a $4$-cycle, we obtain $V(G_v^*) = \{v, w, w', z\}$ and $wz, w'z \subseteq E(G_v^*)$ for some $z \in V(G^*) \setminus \{v, w, w'\}$. Suppose $z \neq y_3$. Since $\Delta(G) = 4$, we have $N(y_1) = \{x_1, x_1^*, y_2, y_4\}$, so $(\{y_3\}, \emptyset)$ is a component of $G^*$. Thus, (\ref{h'=1a}) becomes
\[\iota(G,C_4) \leq \frac{|N[y_1] \cup \{y_3\}| - 1}{5} + \frac{|V(G_v^*)| + 1}{5} + \sum_{H \in \mathcal{I}}\frac{|V(H)|}{5} = \frac{n}{5}.\]
Now suppose $z = y_3$. Let $Z = \{y_3, w, w', y_2, y_4\}$. Let $G_Z = G-Z$. We have $N_{G_Z}(v) = \{x_1, x_1^*\}$ and $y_1 \in N(x_1)$, so $G_Z$ is connected. Since $\delta(G_Z) < 3$, $G_Z$ is not a $\cG_9$-graph. We have $|V(G_Z)| = n - 5 \geq 5$, so $G_Z$ is not an $\cE_4$-graph. By the induction hypothesis, $\iota(G_Z,C_4) \leq |V(G_Z)|/5$. By Lemma~\ref{LemIsol1}, $\iota(G,C_4) \leq 1 + \iota(G_Z,C_4)$ as $Z \subseteq N[y_3]$. Thus, $\iota(G,C_4) \leq |Z|/5 + |V(G_Z)|/5 = n/5$.

Now suppose $|N(y) \cap N(v)| \leq 1$ for each $y \in V(H_1)$. Thus, $y_1\neq y_1'$. Let $G^* = G - N[y_1]$. Then, $G^*$ has a component $G_v^*$ such that $N[v]\setminus \{x_1\} \subseteq V(G_v^*)$. Since $|N(y_3) \cap N(v)| \leq 1$, if $y_3 \in V(G_v^*)$, then $y_3$ is a leaf of $G_v^*$. Let $G_v^- = G_v^* - (\{y_3\} \cap V(G_v^*))$. By Lemma~\ref{LemIsol1}, Lemma~\ref{LemIsol2} and Corollary~\ref{LemIsolG'cor}, 
\begin{align*} \iota(G,C_4) \leq 1 + \iota(G_v^-,C_4) + \sum_{H \in \cH_{x_1}^*} \iota(H,C_4) \leq \frac{|N[y_1]\cup\{y_3\}|}{5} + \iota(G_v^-,C_4) +  \sum_{H\in\cH_{x_1}^*}\frac{|V(H)|}{5}.
\end{align*}
This gives $\iota(G,C_4) \leq n/5$ if $\iota(G_v^-,C_4) \leq |V(G_v^-)|/5$. Suppose $\iota(G_v^-,C_4) > |V(G_v^-)|/5$. By the induction hypothesis, $G_v^-$ is an $\cE_4$-graph. We have $d_{G_v^-}(x_1') \leq |N(x_1') \setminus \{y_1'\}| \leq \Delta(G) - 1 = 3$. If $d_{G_v^-}(x_1') < 3$, then by Observation~\ref{ObsG9}, $G_v^-$ is not a $\cG_9$-graph. Again by Observation~\ref{ObsG9}, if $d_{G_v^-}(x_1') = 3$, then, since $vx_1' \in E(G_v^-)$ and $d_{G_v^-}(v) = 3$, $G_v^-$ is not a $\cG_9$-graph. Thus, $|V(G_v^-)| = 4$, $V(G_v^-) = N[v] \setminus \{x_1\}$, and hence $\cH_{x_1'}^* = \emptyset$. Let $F = G - N[y_1']$ and $F^- = G - (\{x_1'\} \cup V(H_1))$. Then, $F^-$ is a connected graph and $|V(F^-)| = n - 5 \geq 5$. Similarly to $G_v^-$, $F^-$ is not a $\cG_9$-graph, so $F^-$ is not an $\cE_4$-graph. By the induction hypothesis, $\iota(F^-,C_4) \leq |V(F^-)|/5$. By Lemma~\ref{LemIsol1}, $\iota(G,C_4) \leq 1 + \iota(F,C_4) \leq |\{x_1'\} \cup V(H_1)|/5 + \iota(F,C_4)$. Thus, $\iota(G,C_4) \leq n/5$ if $F = F^-$. Suppose $F \neq F^-$. Then, $F^- = F - y$ for some $y \in V(H_1)$. Since $|N(y) \cap N(v)| \leq 1$, $\iota(F,C_4) = \iota(F^-,C_4)$ by Lemma~\ref{LemIsolG'}. Thus, we again have $\iota(G,C_4) \leq |\{x_1'\} \cup V(H_1)|/5 + |V(F^-)|/5 = n/5$.
\qed

\section{Problems}

Theorem~\ref{ThmC4isol} immediately raises certain questions. The second part of the theorem tells us that the bound in the theorem is attained by $B_{n, C_4}$ for every $n$. It is easy to see that many other graphs attain the bound, but it would be interesting to know what structures extremal graphs have. 
This leads us to our first two problems.

\begin{problem} \label{problem1}
Which are the connected graphs that attain the bound in Theorem~\ref{ThmC4isol}?
\end{problem}

\noindent
It is quite likely that Problem~\ref{problem1} is significantly challenging for the case where $n$ is not divisible by $5$. On the other hand, the problem for the case where $n$ is divisible by $5$ should be tractable, and this is the following relaxation of Problem~\ref{problem1}.

\begin{problem}
Which are the connected graphs $G$ such that $\iota(G,C_4) = |V(G)|/5$?
\end{problem}

Let $\mathcal{C}_4$ denote the set $\{C_k \colon k \geq 4\}$. In view of Theorems~\ref{Borgcycle} and \ref{ThmC4isol}, we also pose the following problem. 

\begin{problem} \label{problem3}
What is the smallest real number $c$ (if it exists) such that for some finite set $\mathcal{E}$ of graphs, $\iota(G, \mathcal{C}_4) \leq c |V(G)|$ for every connected graph $G$ that is not an $\mathcal{E}$-graph?
\end{problem}

\noindent
We claim that Theorem~\ref{ThmC4isol} extends as follows.

\begin{conjecture} The number $c$ in Problem~\ref{problem3} is $1/5$.
\end{conjecture}

\appendix
\section*{Appendix: Mathematica code for Lemma~\ref{ThmC4isolsmall}} 

The following Mathematica code outputs the six graphs in $\cG_9$ and makes use of the IGraph/M package \cite{H-IGM} and Brendan McKay's graph data \cite{McKay}.\\

\begin{lstlisting}[language=Mathematica]
MinMaxDegQ[d_, D_, g_] := 
  TrueQ[Min[VertexDegree[g]] >= d && Max[VertexDegree[g]] <= D];
Min2Max4DegQ[g_] := MinMaxDegQ[2, 4, g];
Needs["IGraphM`"];
graph9c = Import["graph9c.g6"];
graph9special = Select[graph9c, Min2Max4DegQ];
F = CycleGraph[4];
Exc = {};
For[j = 1, j <= Length[graph9special], j++,
  G = graph9special[[j]];
  ub = Floor[Length[VertexList[G]]/5];
  S = Subsets[VertexList[G], {ub}];
  nisol = True;
  i = 1;
  While[nisol && i < (Length[S] + 1),
   GNS = Subgraph[G, 
     Complement[VertexList[G], 
      VertexList[NeighborhoodGraph[G, S[[i]], 1]]]];
   If[IGSubisomorphicQ[F, GNS],
    i++,
    nisol = False
    ];
   ];
  If[nisol, Exc = Append[Exc, G], ## &[]];
  ];
Exc
\end{lstlisting}

\vskip 5mm

\noindent
\textbf{Acknowledgements.} The authors are grateful to the two anonymous referees for checking the paper carefully and providing constructive remarks.


\end{document}